\documentclass[11pt,a4paper]{article}
\usepackage{amsmath,amsthm}
\usepackage{amsfonts,amstext,amssymb, mathtools, comment, graphicx}

\newcommand{\bs}{\boldsymbol}
\newcommand{\one}{{\bs 1}}
\newcommand{\p}{{\mathbb P}}
\newcommand{\e}{{\mathbb E}}
\renewcommand{\a}{\alpha}
\newcommand{\D}{\mathrm d}
\newcommand{\levy}{L\'{e}vy }

\newcommand{\R}{\mathbb R}
\renewcommand{\Re}{{\rm Re}}

\newcommand{\1}[1]{\mbox{\large  1}_{\{#1\}}}
\newcommand{\pii}{{\bs \pi}}
\newcommand{\matI}{\mathbb{I}}
\newcommand{\matO}{\mathbb{O}}
\newcommand{\Lmb}{\Lambda}

\newcommand{\bh}{{\bs h}}
\newcommand{\zero}{{\bs 0}}

\newtheorem{thm}{Theorem}
\newtheorem{lem}[thm]{Lemma}
\newtheorem{prop}[thm]{Proposition}
\newtheorem{cor}[thm]{Corollary}

\theoremstyle{remark}
\newtheorem{rem}{Remark}[section]

\begin{document}
\bibliographystyle{plain}
\title{Markov-modulated Brownian motion with two reflecting barriers}
\author{Jevgenijs Ivanovs\footnote{Eurandom, Eindhoven
University of Technology and Korteweg-de Vries Institute for
Mathematics, University of Amsterdam}} \maketitle \begin{abstract}
We consider a Markov-modulated Brownian motion reflected to stay
in a strip $[0,B]$. The stationary distribution of this process is
known to have a simple form under some assumptions. We provide a
short probabilistic argument leading to this result and explaining
its simplicity. Moreover, this argument allows for generalizations
including the distribution of the reflected process at an
independent exponentially distributed epoch. Our second
contribution concerns transient behavior of the reflected system.
We identify the joint law of the processes $t,X(t),J(t)$ at
inverse local times.
\end{abstract}

\vspace{0.1in}\noindent \hspace{-0.1in}\begin{tabular}{l l}
AMS 2000 Subject classification:&Primary 60J55\\&Secondary 60K25;60K37\end{tabular}\\

\section{Introduction}
In this paper we investigate a Markov-modulated Brownian motion
(MMBM) reflected to stay in a strip $[0,B]$, where $B>0$. Roughly
speaking, an MMBM is a process with piecewise Brownian paths with
drift and variance parameters determined by a finite-state Markov
chain. The variance parameters are allowed to be 0, in which case
the corresponding pieces are linear. Mathematically, MMBM is just
a Markov additive process with continuous
paths~\cite[Ch.~XI]{asmussen:apq}. Letting $X(t)$ be an MMBM we
construct a doubly reflected process $W(t)$ as
\begin{equation}\label{eq:reflection_relation}W(t)=X(t)+L(t)-U(t),\end{equation} where
$L(t)$ and $U(t)$ are the \emph{local times} at respectively the lower
and the upper barriers (that is at 0 and at $B$), given as the
solutions of a Skorokhod problem, see Section~\ref{sec:reflection}
for details.

The model considered in this paper is also called in the
literature \emph{a second-order fluid model} or \emph{a fluid
model with Brownian noise}. It was introduced as a generalization
of an extensively studied fluid flow model, where it is assumed
that all the variance parameters are 0, making the process
piecewise linear. Fluid models were initially proposed for
manufacturing and telecommunication systems, where units of work
(products or packets) are processed so fast that they can be
modelled as fluid instead of discrete units. Since then the use of
fluid models has become widespread making it a classical model in
applied probability with a variety of application areas like the
theory of queues and dams, risk processes and insurance,
environmental problems, etc. The literature on this topic is
extensive, we only mention the seminal
papers~\cite{kosten,anick_mitra_sondhi}, a
survey~\cite{fluid_survey}, and a more recent
paper~\cite{ahn_badescu_ramaswami} with an extended list of
references.

The second-order fluid models were simultaneously introduced
in~\cite{asmussen_fluid_flow,karandikar_kulkarni,rogers}. The
highly cited paper~\cite{rogers} by L.~C.~G.~Rogers can be
considered as one of the most influential papers, not just in the
theory of fluid models, but in the much more general theory of
fluctuations of Markov additive processes. In~\cite{rogers} the
stationary distribution of a reflected MMBM is derived for both a
single barrier and two barriers assuming that either all the
variance parameters are 0 or all are positive, see also the comments
in Section~\ref{sec:rogers}. The case of a single barrier, see
also~\cite{asmussen_fluid_flow,karandikar_kulkarni}, is a special
case of the latter problem with $B=\infty$ and is substantially
easier to analyze. In fact, in this case the analysis can often
be extended to more general Markov-modulated \levy models with
jumps, where all the jumps are either positive or negative; see for example~\cite{prabhu_zhu}, where the stationary distribution
of an infinite buffer Markov-modulated M/G/1 queue is obtained.

The form of the stationary distribution of an MMBM with two
reflecting barriers is strikingly simple. It is noted
in~\cite{rogers}, where the result is derived using the theory of
generators of Markov processes, that there has to be a
probabilistic story explaining why the solutions to the
differential equations reduce so far. We provide a very simple
probabilistic argument, which offers this explanation. Moreover,
this argument allows us to find the stationary distribution of a
reflected MMBM with both linear and Brownian segments, which to
our knowledge has not been done so far (its transform is obtained
in~\cite{lambda}). Furthermore, we obtain the distribution of the
reflected process at an independent exponentially distributed
epoch, which gives a first insight in the transient behavior of
the reflected process.

Our second contribution concerns the transient behavior of the
reflected system. We identify the joint law of the processes
$t,X(t),J(t)$ at inverse local times. We envisage that this key
result will find various applications; some are mentioned in
Section~\ref{sec:localtimes}. In particular, we show that a result
from~\cite{lambda} on stationary overflow and unused capacity is a
simple consequence of our result. Moreover, considering the
special case of a Brownian motion we recover the results
from~\cite{williams_RBM}.

Fluid models play a prominent role in applied probability. The
importance comes from the fact that they are flexible enough to
model a variety of different phenomena, and at the same time the
analysis often remains tractable. One important observation is
that phase-type jumps can be easily incorporated in the model.
This is done by introducing auxiliary states of $J(t)$ and
`levelling' the jumps out, see e.g.~\cite{pistorius}. Notice that
this observation shows the importance of an MMBM, where some
variance parameters are allowed to be 0 (the case not analyzed so
far in the literature). We also remark that the level $B$ of the
barrier can be interpreted in different ways. In the queueing
context it is the buffer size. In the financial setting it can,
e.g., represent the height of a dividend barrier, in which case
$U(t)$ represents the cumulative amount of dividends paid before
time $t$. Throughout this paper we use the queueing terminology.

This paper has two parts: Section~\ref{sec:statdistr} containing
the results about the stationary distribution, and
Section~\ref{sec:localtimes} discussing the inverse local times.
Some basic facts and definitions concerning two-sided reflection
and MMBM are given in Section~\ref{sec:prelim}.

\section{Preliminaries}\label{sec:prelim}
Throughout this work we use bold symbols to denote vectors
(usually column-vectors), in particular $\one$ and $\zero$ are
vectors of 1s and 0s respectively. We write $\matI$ for an
identity matrix and $\Delta_{\bs v}$ for a diagonal matrix with
entries of $\bs v$ on the diagonal.
\subsection{Two-sided reflection}\label{sec:reflection}
Let $X(t), t\geq 0$ be a real continuous function with $X(0)\in
[0,B]$ (a fixed sample path of a stochastic process). The
two-sided reflection $W(t)$ of $X(t)$, with respect to the strip
$[0,B]$, is defined through (\ref{eq:reflection_relation}), where
$W(t),L(t),U(t)$ are real continuous functions which satisfy the
following conditions:
\begin{itemize}
  \item $L(t)$ and $U(t)$ are non-decreasing with $L(0)=U(0)=0$,
  \item $W(t)\in[0,B]$ for all $t\geq 0$,
  \item $W(s)=0$ if $\forall t>s: L(s)<L(t)$, and $W(s)=B$ if $\forall
  t>s:
  U(s)<U(t)$.
\end{itemize}
The last condition states that the points of increase of $L$ and
$U$ are contained in $\{t\geq 0:W(t)=0\}$ and $\{t\geq 0:W(t)=B\}$
respectively. It is known that such a triplet of functions exists
and is unique~\cite{ramanan_skorohod0a}, and is called the
solution of the two-sided Skorokhod problem. The functions $L(t)$
and $U(t)$ are called local times at the lower and upper barriers
respectively, that is at $0$ and at $B$.

\subsection{MMBM and its basic properties}
We start with a rigorous definition of an MMBM, see also
~\cite[Ch.~XI]{asmussen:apq}. Let $J(t)$ be a continuous-time
irreducible Markov chain with a finite state space $E$, where
$|E|=N$, and transition rate matrix $Q$. Denote its invariant
distribution through $\pii$. For each $i\in E$ let $X_i(t)$ be a
Brownian motion with variance $\sigma_{i}^{2}\geq 0$ and drift
$\mu_{i}$; the vectors with elements $\sigma_i$ and $\mu_i$ are
denoted through $\bs \sigma$ and $\bs \mu$ respectively. We assume
that $J(t)$ and all $X_{i}(t)$ are mutually independent. Letting
$T$ be the last jump epoch of $J$ before $t$ ($T=0$ if there were
no jumps), we define $X(t)$ recursively as
\begin{equation}\label{eq:mmbm_def}X(t)=X(T)+\sum_{i\in E}\1{J(T)=i}(X_{i}(t)-X_{i}(T))\text{ and }X(0)=0.\end{equation}
Informally, $X(t)$ evolves as $X_i(t)$ while $J(t)$ is in
state $i$. The law of $(X(t),J(t))$ given $J(0)=i$ is denoted
through $\p_i$, and the expectation through $\e_i$. Sometimes we
consider the process $(x_0+X(t),J(t))$ for an arbitrary $x_0\in\R$
and still write $(X(t),J(t))$ where no confusion can arise.

Let us now introduce some important concepts related to MMBM.
The asymptotic drift of the process $X(t)$ is denoted through
$\kappa:=\lim_{t\rightarrow\infty}X(t)/t=\sum_{i\in E}\pi_{i}\mu_{i},$ which does not depend on the initial state of $J(t)$.
Next, define the subsets $E^{+}$ and $E^{-}$ of $E$ through
\[E^{\pm}=E\backslash\{i\in E:\sigma^{2}_{i}=0,\pm\mu_{i}\leq 0\}.\] Thus, e.g., $E^{+}$ is the set of indices
of the processes $X_{i}(t)$ which are not monotonically decreasing; $X_i(t),i\in E^+$ are the processes which hit the set $(0,\infty)$ with positive probability. In addition, we let
$N^\pm=|E^{\pm}|$ and write $M_\pm$ to denote a restriction of an arbitrary matrix $M$ with $N$ rows obtained by keeping those and only those rows which are indexed by $E^\pm$.

Let $e_{q}$ be an exponential random variable with rate $q\geq 0$,
independent of everything else, where by convention
$e_{0}=\infty$. It is often convenient to consider the process
$(X(t),J(t))$ only until $e_q$, at which moment it is sent to some
additional absorbing state $(\partial_X,\partial_J)$, in other
words - killed. The law of such a process is denoted through
$\p^q$. In order to simplify notation, we also keep writing simply
$\p$ if there is no confusion. We do so in the rest of this
section.

Define the first passage times for $x\geq 0$ as
\begin{equation}\nonumber
\tau_x^\pm=\inf\{t\geq 0:\pm X(t)>x\}.
\end{equation}
Note that $\tau_x^+=\infty$ if $X(t)$ does not exceed level $x$
before $e_q$. The continuity of paths and the strong Markov
property of MMBM, and in addition the memoryless property of the
exponential distribution, imply that the time changed process
$J(\tau_x^+)$ is a Markov chain itself, see also~\cite{rogers}.
Note that $J(\tau^+_x)$ jumps to the absorbing state $\partial_J$
at $x=\sup\{X(t):0\leq t<e_q\}$. Apart from $\partial_J$ this
Markov chain can take values only in $E^+$. Furthermore it is
irreducible when restricted to the set $E^+$. Let $\Lmb^+$ be the
corresponding $N^+\times N^+$ transition rate matrix of
$J(\tau^+_x)$. Let also $\Pi^+$ be a $N\times N^+$ matrix with elements
$\p_i(J(\tau^+_0)=j)$, where $i\in E,j\in E^+$ (with an obvious order); we often simply write $\p(J(\tau_0^+))$.
Using this notation we can furthermore write \[\p(J(\tau_x^+))=\Pi^+e^{\Lmb^+x}.\]
The Markov chain $J(\tau_{x}^{+})$ (restricted to $E^+$) is recurrent if and only if $q=0$ and $\kappa\geq 0$.
Finally, we consider $\tau^-_x$ and define
the matrices $\Lmb^-$ and $\Pi^-$ in an obvious way.
Note that the matrices $\Lambda^\pm$ and $\Pi^\pm$ depend on $q$. In order to stress this dependence and distinguish from the special case of $q=0$, we sometimes call them $q$-killed versions.

The matrices $\Lambda^\pm$ and $\Pi^\pm$ play a prominent role in
the following. It is, therefore, required in practice to be able
to compute these matrices. There are two methods in the
literature. Firstly, one can use an iteration scheme, see
e.g.~\cite{asmussen_fluid_flow,breuer}. Secondly, a Jordan normal
form of $\Lambda^\pm$ can be obtained directly using the theory of
generalized Jordan chains, see~\cite{lambda,time_rev}. Finally, it
is known that
\begin{equation}\label{eq:lambda}\frac{1}{2}\Delta^2_{\bs
\sigma}\Pi^\pm(\Lambda^\pm)^2\mp\Delta_{\bs
\mu}\Pi^\pm\Lambda^\pm+(Q-q\matI)\Pi^\pm=\matO,\end{equation}
where $\matO$ is $N\times N^\pm$ matrix of zeros, see for
example~\cite{breuer}.

\section{The stationary distribution}\label{sec:statdistr}
In the first part of this section we provide an elementary
derivation of the stationary distribution of $(W(t),J(t))$; we let
$(W,J)$ denote the corresponding random vector. Let $\hat J(t)$ be the
time-reversed version of $J(t)$, that is, $\hat J(t)$ is a Markov
chain with transition rate matrix
$\Delta_{\pii}^{-1}Q^T\Delta_{\pii}$. Define $\hat X(t)$ according
to (\ref{eq:mmbm_def}) with $\hat J(t)$ in place of $J(t)$. The
law of the time-reversed MMBM is denoted through $\hat\p$. We
start with the following basic identity
\begin{equation}\label{eq:stat_dist00}\p(W\geq x|J=i)=\hat\p_i(X(\tau_{[x-B,x)})\geq x),\end{equation}
where $\tau[u,v)=\inf\{t\geq 0:X(t)\notin [u,v)\}$. This type of
representation was first noted in~\cite{lindley59} in the case of
a random walk with two reflecting barriers. A short derivation of
its continuous-time analogue is given in~\cite[Prop.~3.7,
Ch.~XIV]{asmussen:apq}, see also~\cite{loss_rates} for the case of
Markov additive input. It is well-known that an MMBM can not hit a
level without passing it, so we obtain
\begin{equation}\label{eq:stat_dist_0}\hat\p(W\geq x|J=i)=\p_i(\tau_x^+<\tau_{B-x}^-),\text{ where }x\in (0,B].\end{equation}
Note that this identity indeed does not hold for $x=0$ and $i\in E\backslash E^+$.

Next, for $a,b\geq 0$ with at least one being strictly positive
consider the matrices
\begin{align}\nonumber &C(a,b)=\p(\tau^+_{a}<\tau_b^-,J(\tau_{a}^+)) &D(a,b)=\p(\tau_{b}^-<\tau_a^+,J(\tau_{b}^-))\end{align} of
dimensions $N\times N^+$ and $N\times N^-$ respectively.
The strong Markov property implies
\begin{align*}C(a,b)&=\p (J(\tau_a^+))-\p
(\tau_b^-<\tau_a^+,J(\tau_a^+))=\p
(J(\tau_a^+))-D(a,b)\p_-(J(\tau_{a+b}^+)),\end{align*}
where $\p_-$ denotes a restriction according to $J(0)\in E^-$.
Expand $D(a,b)$ in a similar way to arrive at
\begin{align}\label{eq:C_D}
&C(a,b)=\Pi^+e^{a\Lambda^+}-D(a,b)\Pi^+_-e^{(a+b)\Lambda^+},\\
\nonumber&D(a,b)=\Pi^-e^{b\Lambda^-}-C(a,b)\Pi^-_+e^{(a+b)\Lambda^-}.
\end{align}
Assuming that $\kappa\neq 0$ we denote
\begin{align}\label{eq:K_def}&K^+:=(\matI-\Pi^-_+e^{B\Lambda^-}\Pi^+_-e^{B\Lambda^+})^{-1}&K^-:=(\matI-\Pi^+_-e^{B\Lambda^+}\Pi^-_+e^{B\Lambda^-})^{-1},\end{align}
where the inverses are well defined, because one of the matrices
$\Pi^-_+e^{B\Lambda^-}$ and $\Pi^+_-e^{B\Lambda^+}$ is a transient
probability matrix (and the other is recurrent). Thus
\begin{equation}\label{eq:C}C(x,B-x)=(\Pi^+e^{x\Lambda^+}-\Pi^-e^{(B-x)\Lambda^-}\Pi^+_-e^{B\Lambda^+})K^+\end{equation}
and then (\ref{eq:stat_dist_0}) leads to the following result
\begin{thm}
If $\kappa\neq 0$ then for $x\in (0,B]$ it holds that
\begin{equation}\label{eq:stat_dist}\hat\p(W\geq x|J)^T=(\Pi^+e^{x\Lambda^+}-\Pi^-e^{(B-x)\Lambda^-}\Pi^+_-e^{B\Lambda^+})K^+\one.\end{equation}
\end{thm}
Finally, we give a brief comment on the more delicate case of $\kappa=0$.
The problem of this case lies in the fact that the system of equations (\ref{eq:C_D}) does not identify the matrices $C(a,b)$ and $D(a,b)$ uniquely. In fact, exactly one equation is missing.
This equation is of the form
\[C(a,b)(a\one_++\bh_+)+D(a,b)(-b\one_-+\bh_-)=\bh,\]
where $\bh$ is any vector which solves $Q\bh+{\bs \mu}=\zero$,
see~\cite[Section 7]{lambda}.

\subsection{Further identities}
First, we present a very simple coupling argument, which allows us
to derive an equation that complements~(\ref{eq:stat_dist}).
Consider the two-sided reflection $\tilde W(t)$ of $(-X(t),J(t))$
in $[0,B]$. Assuming $X(0)=0$ and $\tilde X(0)=B$ it is easy to
see that $\p(W(t)\leq x|J(t))=\p(\tilde W(t)\geq B-x|J(t))$.
Letting $t\rightarrow\infty$ we obtain from~(\ref{eq:stat_dist})
\begin{equation}\label{eq:stat_dist1}\hat\p(W\leq x|J)^T=\hat\p(\tilde W\geq B-x|J)^T=(\Pi^-e^{(B-x)\Lambda^-}-\Pi^+e^{x\Lambda^+}\Pi^-_+e^{B\Lambda^-})K^-\one,\end{equation}
where $x\in [0,B)$. Note that $\Lambda^+,\Pi^+$ become $\Lambda^-,\Pi^-$, because the MMBM $(-X(t),J(t))$ is used to construct $(\tilde W,J)$.

Note that the equations (\ref{eq:stat_dist}) and (\ref{eq:stat_dist1})
lead to two different representations of the density $\hat \p(W\in
\D x|J)$. In addition, one easily obtains the point masses at $0$
and $B$ (from~(\ref{eq:stat_dist1}) and~(\ref{eq:stat_dist}) respectively) :
\begin{align}\nonumber \label{eq:point_masses}&\hat\p(W=0|J)^T=(\Pi^--\Pi^+\Pi^-_+)e^{B\Lambda^-}K^-\one, &\hat\p(W=B|J)^T=(\Pi^+-\Pi^-\Pi^+_-)e^{B\Lambda^+}K^+\one.\end{align}
Finally, $\Pi_+^+=\matI$ and $\Pi^-_-=\matI$ and hence $W$ has no mass at 0 (respectively $B$) given $J$ is in $E^+$ (respectively $E^-$).

\subsection{On the result of Rogers~\cite{rogers}}\label{sec:rogers}
In this section we comment on the result of Rogers~\cite{rogers}.
We are only interested in a fluid model with Brownian noise and a
finite buffer presented in~\cite[Section 7]{rogers}. It is assumed
there that the fluid evolves as an independent sum of a Markov
modulated linear drift and a standard Brownian motion. At first
sight, this is a rather special case of an MMBM. Note, however,
that the process $(X(t),J(t))$ can be time-changed without
changing the distribution of $\left(W|J\right)$ in the following
way. We scale time by $c_i>0$ while $J(t)$ is in state $i$, that
is, we consider a new MMBM specified by the transition rate matrix
$\Delta_{\bs c}^{-1}Q$ and parameters $\sigma^2_i/c_i,\mu_i/c_i$.
It is easy to see that this new MMBM gives rise to the same
distribution of $\left(W|J\right)$. Hence,~\cite{rogers}, in fact,
does not assume more than this: all the variance parameters are
strictly positive. In other words $E^+=E^-=E$ and hence
$\Pi^\pm=\matI$. Under this assumption (\ref{eq:stat_dist})
results in the following:
\[\p(W\in\D
x|J)^T=-(e^{x\hat\Lambda^+}\hat\Lambda^++e^{(B-x)\hat\Lambda^-}\hat\Lambda^-e^{B\hat\Lambda^+})(\matI-e^{B\hat\Lambda^-}e^{B\hat\Lambda^+})^{-1}\one,\]
which is (7.13) of~\cite{rogers} up to the minus sign. Here
$\hat\Lambda^\pm=\hat\Gamma_\mp$, because of the different
definition of time-reversal (3.3) of~\cite{rogers}. The missing
minus sign is a consequence of a mistake in normalization
in~\cite{rogers}. Finally, we note that in~\cite{rogers} the above
result was obtained using the theory of generators of Markov
processes. It required solving second-order differential
equations, and verification of the solution (positivity of the
density is not established in~\cite{rogers} though). On the
contrary, our result was obtained directly using elementary
probabilistic arguments. This allowed us to treat the problem in
its full generality, i.e., no assumption on variances.

\subsection{The distribution at an exponential epoch}
In this section we identify $\p_i(W(e_q),J(e_q)=j)$ for
$X(0)\in\{0,B\}$ (start at a boundary), which provides some
information of the transient behavior of the reflected process. An
inspection of the proof of~\cite[Prop.~3.7, Ch.~XIV]{asmussen:apq}
reveals that a representation similar to~(\ref{eq:stat_dist00})
holds true for finite time $T$:
\[\p_i(W(T)\geq x|J(T)=j)=\hat\p_j(\tau_{[x-B,x)}\leq T,X(\tau_{[x-B,x)})\geq x|J(T)=i),\]
where it is assumed that $X(0)=0$. Note that
$\pi_i\p_i(J(T)=j)=\pi_j\hat\p_j(J(T)=i)$ to arrive at the
following equation
\[\hat\p_i(W(e_q)\geq x,J(e_q)=j)=\p_j(\tau_x^+<e_q,\tau_x^+<\tau_{B-x}^-,J(e_q)=i)\frac{\pi_j}{\pi_i},\]
which in matrix form can be written as
\[\hat\p(W(e_q)\geq x,J(e_q))^T=\Delta_\pi\p^q(\tau_x^+<\tau_{B-x}^-,J(\tau_x^+))\p(J(e_q))\Delta_\pi^{-1}.\]
Moreover, one can show that (\ref{eq:C}) holds true under exponential killing.
Noting that $\p(J(e_q))=q(q\matI-Q)^{-1}$ we find
\[\hat\p(W(e_q)\geq x,J(e_q))^T=q\Delta_\pii(\Pi^+e^{x\Lambda^+}-\Pi^-e^{(B-x)\Lambda^-}\Pi^+_-e^{B\Lambda^+})K^+(q\matI-Q)^{-1}\Delta_\pii^{-1},\]
where $x\in(0,B]$ and all the occurrences of matrices
$\Lambda^\pm$ and $\Pi^\pm$ refer to the $q$-killed versions.
Finally, one can derive a symmetric equation for the case
$X(0)=B$. The distribution of $W(e_q)$ in the case of
$X(0)\in(0,B)$ does not have an explicit form. In this case one
needs to resort to its Laplace transform as in~\cite{lambda}.

\section{Inverse local times}\label{sec:localtimes}
This section is devoted to the study of the transient behavior of
the reflection system\\ $(t,X(t),J(t),W(t),L(t),U(t))$. More
concretely, we characterize the joint law of these processes at
inverse local times $\tau^L_x,x\geq 0$ and $\tau^U_x,x\geq 0$,
where
\[\tau^L_x=\inf\{t\geq 0:L(t)>x\}\text{ and }\tau^U_x=\inf\{t\geq 0:U(t)>x\}\]
for any $x\geq 0$. This key result allows us to answer a number of
important questions. For example, given $X(0)\in [0,B]$ and
$J(0)\in E$, when does the buffer become empty for the first time and what is the
state of $J(t)$ at this time? What is the amount of lost fluid
until then? Mathematically speaking, we are interested in
$(\tau_0^L,J(\tau_0^L),U(\tau_0^L))$. We can also ask: what is the
length of an arbitrary busy period? What is the amount of lost
fluid during a busy period given there was a loss? Moreover, we
condition on the state $i$ of $J(t)$ right before this busy period
starts and the state $j$ at which it finishes. These quantities
are described by the jumps of $\tau_x^L$ and $U(\tau_x^L)$ given
there is a corresponding transition of $J(\tau_x^L)$. The answer
to these questions is immediate in view of
Theorem~\ref{thm_localtime}. Moreover, in
Section~\ref{sec:overflow} we show that the stationary overflow and
unused capacity vectors can be trivially obtained from
Theorem~\ref{thm_localtime}; by doing so we
recover a result from~\cite{lambda}. Finally, in
Section~\ref{sec:special_cases} we consider a special case of a
simple Brownian motion and recover the results
from~\cite{williams_RBM}, where as an easy consequence an
asymptotic variance of the overflow process is obtained.

It should be noted that the analysis of the two-sided reflection
problem, compared to the one-sided one, is considerably harder.
The main problem lies in the fact that there are \emph{two} local
time processes, which are interrelated in an intricate way. The
crucial idea is to study the set of points $x\geq 0$
 such that $X(\tau_x^L)=y$ for a fixed $y\in\R$, see Lemma~\ref{lem:tau_sigma}.
 This idea in a simpler form also appears in~\cite[Section 5]{rogers}, where it was used to derive point masses of $W(t)$ at $0$ and $B$ in stationarity.

\subsection{Markov additive processes}\label{sec:MAP}
It will be shown that $(X(t),J(t))$ observed at inverse local
times is a Markov additive process (MAP). In this section we
present a definition and some basic properties of a MAP. A MAP is
a bivariate Markov process $(Y(t),\tilde J(t))$, which satisfies
the following property for every $T\geq 0$. Given $\tilde J(T)=i$
the shifted process $(Y(T+t)-Y(T),\tilde J(T+t)),t\geq 0$ is
independent from $(Y(t),\tilde J(t)),0\leq t\leq T$ and has the
same law as $(Y(t),\tilde J(t)),t\geq 0$ given $\tilde J(0)=i$. It
is commonly assumed that the state space of $\tilde J(t)$ is
finite, in our case it is a subset of $E$. Such a process can be
seen as a Markov-modulated \levy process with additional jumps at
the transition epochs of $\tilde J(t)$,
see~\cite[Ch.~XI]{asmussen:apq}. Hence a MAP with continuous paths
is an MMBM. Importantly, if $Y(t)$ has no negative jumps, then for
all $\a\leq 0$ there exists a square matrix function $F(\a)$ such
that
\begin{equation}\label{eq:ij}\e[e^{\a Y(t)};\tilde
J(t)]=e^{F(\a)t},\end{equation} where the term on the left side denotes a matrix with elements $\e_i[e^{\a Y(t)}\1{\tilde
J(t)=j}]$ for $i,j$ in the state space of $\tilde J(t)$. One can see now that $F(\a)$ uniquely specifies the
law of the process $(Y(t),\tilde J(t))$. This matrix can be written explicitly
in terms of the \levy exponents of the underlying \levy processes,
the Laplace transforms of the additional jumps (at transition
epochs), and the transition rate matrix of $\tilde J(t)$. Finally,
there exists a (Perron-Frobenius) eigenvalue $k(\a)$ of $F(\a)$,
which is real and is larger than the real parts of all the other
eigenvalues of $F(\a)$.

\subsection{The main result}
We start by making the following observations:
\begin{itemize}
  \item $\tau^L_x$ and hence also $J(\tau^L_x), X(\tau^L_x), U(\tau_x^L)$ are right-continuous;
  \item $L(\tau_x^L)=x$ by the continuity of $L(t)$;
  \item $W(\tau_x^L)=0$, because $\tau_x^L$ is a point of increase of $L(t)$;
  \item $U(\tau_x^L)$ is piece-wise constant;
  \item $X(\tau_x^L)=U(\tau_x^L)-x$ is piece-wise linear.
\end{itemize}
Moreover, by the strong Markov property of $(X(t),J(t))$ it follows that $(X(\tau_x^L),J(\tau_x^L))$ is a MAP,
and in particular $J(\tau_x^L)$ is a Markov chain. Note also that $J(\tau_x^L)$ is an irreducible Markov chain taking values in $E^-$.
The additive component $X(\tau_x^L)$ has no negative jumps, hence there exists a $N^-\times N^-$ matrix function $F^L(\a)$ for all $\a\leq 0$, such that
\begin{equation}\label{lt_eq1}\e_- [e^{\a
X(\tau_x^L)};J(\tau_x^L)]=e^{F^L(\a)x},\end{equation}
see Section~\ref{sec:MAP}.
 Clearly,
similar statements hold true with respect to $\tau^U_x$.

Consider the reflected system under the time change $t=\tau_x^L$ (and similarly
$t=\tau_x^U$) and note using the above observations that it is enough to
characterize the trivariate process $(\tau_x^L,X(\tau_x^L),J(\tau_x^L))$.
This process is a MAP with 2-dimensional additive component, and
hence is uniquely specified by the following quantity
\begin{equation}\label{lt_eq2}\e_{x_0}[e^{\a
X(\tau_x^L)-q\tau_x^L};J(\tau_x^L)]=\e_{x_0}[e^{\a
X(\tau_0^L)-q\tau_0^L};J(\tau_0^L)]e^{F^L(\a,q)x},\end{equation}
where $\e_{x_0}[e^{\a X(\tau_0^L)-q\tau_0^L};J(\tau_0^L)]$
describes the initial distribution, and $F^L(\a,q)$ is the
corresponding matrix exponent. Note that $q\geq 0$ can be seen as
the rate of an independent exponential killing of the original
process. The formula (\ref{lt_eq2}) generalizes (\ref{lt_eq1}) by
allowing for arbitrary $q\geq 0$, $X(0)=x_0\in[0,B]$ and $J(0)\in
E$. In the following we do not explicitly write the killing rate
$q$ in order to simplify notation.

As it was shown above, the joint law of $(t,X(t),J(t),L(t),U(t))$ observed at $t=\tau^L_x,x\geq 0$ is uniquely characterized by the quantities $\e_{x_0}[e^{\a
X(\tau_0^L)};J(\tau_0^L)]$ and $F^L(\a)$ ($q$ is implicit here).
Hence our goal is to determine these quantities, as well as those corresponding to $\tau^U_x$.
Define $N\times N^-$ and $N\times N^+$ dimensional matrices $M^L(\a)$ and $M^U(\a)$ through
\begin{align*}&M^L(\a)=\e_{x_0}[e^{\a
X(\tau_0^L)};J(\tau_0^L)](F^L(\a))^{-1},&M^U(\a)=\e_{x_0}[e^{\a
X(\tau_0^U)};J(\tau_0^U)](F^U(\a))^{-1}\end{align*} for those
values of $q\geq 0$ and $\a\leq 0$ for which the inverses are
well-defined. Letting $x_0=0$ and restricting the rows of
$M^L(\a)$ to $E^-$ we get $(F^L(\a))^{-1}$. Hence $\e_{x_0}[e^{\a
X(\tau_0^L)};J(\tau_0^L)]$ and $F^L(\a)$ are readily obtained from
$M^L(\a)$ (a similar statement is true with respect to
$\tau^U_x$). We, therefore, aim to determine the matrices
$M^L(\a)$ and $M^U(\a)$.

Let $\rho^+,\rho^-$ and $k^L(\a),k^U(\a)$ be the Perron-Frobenius
eigenvalues of $\Lambda^+,\Lambda^-$ and $F^L(\a),F^U(\a)$
respectively. It is well-known that $\rho^+,\rho^-\leq 0$. If $q>0$ then the inequalities are strict. If $q=0$ and $\kappa\neq 0$ then one of the inequalities is strict depending on the sign of $\kappa$. In the following we exclude the exceptional case of $q=0$ and $\kappa=0$.
We are ready to present our main result.
\begin{thm}\label{thm_localtime}
Let $\a\in (\rho^-,-\rho^+)$. Then $k^L(\a),k^U(\a)<0$ and
$M^L(\a),M^U(\a)$ are uniquely specified by
\begin{equation}\label{eq:main}[M^L(\a),M^U(\a)]\begin{pmatrix}
  \matI & -\Pi^+_-e^{B\Lambda^+} \\
  -\Pi^-_+e^{B\Lambda^-} & \matI
\end{pmatrix}=[\Pi^-(\Lambda^--\a\matI)^{-1}e^{x_0\Lambda^-},\Pi^+(\Lambda^++\a\matI)^{-1}e^{(B-x_0)\Lambda^+}].\end{equation}
\end{thm}
We make some comments concerning this theorem. Firstly, $k^L(\a)<0$ and $k^U(\a)<0$ ensure that the matrices $M^L(\a)$ and $M^U(\a)$ are well-defined. Secondly, a simple algebraic manipulation shows that~(\ref{eq:main})
can be equivalently rewritten as
\begin{align}\label{eq:thm1}
M^L(\a)=\left(\Pi^-(\Lambda^--\a\matI)^{-1}e^{x_0\Lambda^-}+\Pi^+(\Lambda^++\a\matI)^{-1}e^{(B-x_0)\Lambda^+}\Pi^-_+e^{B\Lambda^-}\right)K^-,\\
\nonumber
M^U(\a)=\left(\Pi^+(\Lambda^++\a\matI)^{-1}e^{(B-x_0)\Lambda^+}+\Pi^-(\Lambda^--\a\matI)^{-1}e^{x_0\Lambda^-}\Pi^+_-e^{B\Lambda^+}\right)K^+,\end{align}
where
$K^-$ and $K^+$ are given in~(\ref{eq:K_def}) and are well-defined unless $q=\kappa=0$, which was excluded from our consideration.

\begin{rem}
Theorem~\ref{thm_localtime} can be rewritten in a very concise way
using a generalized Jordan pair $(V,\Gamma)$ of the analytic
matrix function
$1/2\Delta_{\bs\sigma}^2\a^2+\Delta_{\bs\mu}\a+Q-q\matI$,
see~\cite{lambda}. In particular, using Lemma 6.3
from~\cite{lambda} we arrive at
\[[M^L(\a),M^U(\a)]\begin{pmatrix}
  -V_- \\
  V_+e^{B\Gamma}\end{pmatrix}=[V(\a\matI-\Gamma)^{-1}e^{x_0\Gamma}].\]
This expression may be used in practice when one is interested in
computing the matrices $M^L(\a)$ and $M^U(\a)$.
\end{rem}

\subsection{Proof}
\begin{figure}[h!]
  \centering
  \includegraphics[height=1.5in]{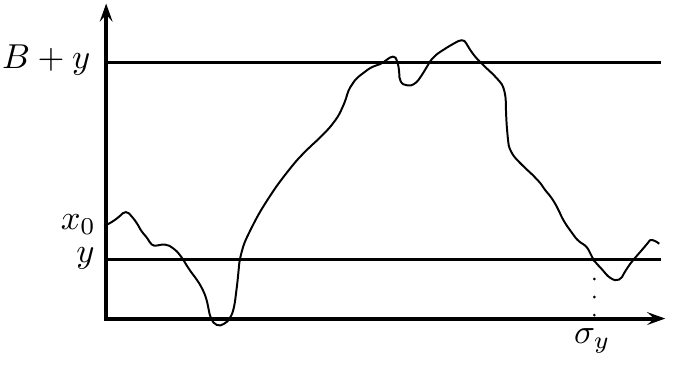}
  \caption{A sample path of $X(t)$.}
  \label{fig:sigma_y}
  \end{figure}

The crucial idea of the proof of Theorem~\ref{thm_localtime} is to
consider the points $x\geq 0$ such that $X(\tau^L_x)=y$ for a
fixed $y\in \R$. Hence we define $x^{(0)}=\inf\{x\geq
0:X(\tau_x^L)=y\}$ and $x^{(n)}=\inf\{x>x^{(n-1)}:X(\tau_x^L)=y\}$
for $n\geq 1$. Recall that $X(\tau_x^L)=U(\tau_x^L)-x$ is
piecewise linear with slope $-1$, so $x^{(n)}$ is strictly larger
than $x^{(n-1)}$. Equivalently, we can look at the time points
$t\geq 0$, such that the local time process $L(\cdot)$ is
increasing and $X(\cdot)$ is at a fixed level $y$ at the time $t$.
The following lemma provides the connection between the above
mentioned points and some quantities which are easily computable.
\begin{lem}\label{lem:tau_sigma}Let \[\varsigma_y=\inf\{t>\tau_{B+y}^+:X(t)<y\}.\]
It holds a.s.~that
\begin{equation}\nonumber\tau_{x^{(0)}}^L=\begin{cases}\tau_{|y|}^-&\text{if }y\leq 0\\
\varsigma_y&\text{if }y>0\end{cases}.\end{equation}
Moreover, for $y=0$ and $X(0)=0,J(0)\in E^-$ it holds a.s.~that
\[\tau_{x^{(1)}}^L=\varsigma_0.\]
\end{lem}
Let us provide some explanation of this result, see also
Figure~\ref{fig:sigma_y}.
If $y\leq 0$ then the first passage time of the level $y$ must be
a point of increase of $L(t)$, otherwise $W(t)$ becomes negative.
If $y>0$ then the first passage time of the level $y$ may not be a
point of increase of $L(t)$. It is necessary that the buffer is
empty when $X(t)$ passes level $y$, which is only possible if an
overflow has occurred before. Hence $X(t)$ should drop by at least
$B$ at the time of hitting $y$, which implies
$\tau_{x^{(0)}}^L=\varsigma_y$. In order to characterize
$\tau_{x^{(n)}}^L,n\geq 1$ we use the strong Markov property.
Hence we only need to consider the case of $y=0$ and
$X(0)=0,J(0)\in E^-$. In this case $\tau_{x^{(0)}}^L=0$ and
$\tau_{x^{(1)}}^L=\varsigma_0$ by a similar argument as above. We
only present a rigorous proof of this latter result.
\begin{proof}[Proof of Lemma~\ref{lem:tau_sigma}]
We assume that $y=X(0)=0,J(0)\in E^-$ and let
$\tau:=\tau^L_{x^{(1)}}$. First, we show that $\tau<\infty$
implies $\varsigma_0\leq \tau$. Observe that
$0=W(\tau)=X(\tau)+L(\tau)-U(\tau)$ and so $U(\tau)=L(\tau)>0$,
because $\tau>\tau_0^-=0$. Hence there was reflection from above
before $\tau$. Let $\hat\tau=\sup\{t\in (0,\tau):W(t)=B\}$ then
$W(\hat\tau)=B,U(\tau)-U(\hat\tau)=0$ and $L(\tau)-L(\hat\tau)\geq
0$, because $L$ is non-decreasing. Thus
$B=W(\hat\tau)-W(\tau)=X(\hat\tau)+L(\hat\tau)-L(\tau)$ and so
$X(\hat\tau)\geq B$. But $X(t)$ can not hit $B$ without passing it
with probability 1, hence $\varsigma_0\leq\tau$ a.s.

 Using the
first part, note that if $\varsigma_0=\infty$ then $\tau=\infty$.
Assuming $\varsigma_0<\infty$ one can easily see that
$\varsigma_0$ is a point of increase of $L(t)$. But
$X(\varsigma_0)=0$ so by the definition of $\tau$ we have
$\tau\leq \varsigma_0$, which concludes the proof.
\end{proof}

Observe that
\begin{align}\nonumber \p(J(\tau^-_{|y|}))&=\Pi^-e^{(x_0-y)\Lambda^-} &\text{ if }y\leq 0,\\ \nonumber\p(J(\varsigma_y))&=\Pi^+e^{(B+y-x_0)\Lambda^+}\Pi^-_+e^{B\Lambda^-} &\text{ if }y\geq 0,
\end{align} where $x_0\in[0,B]$. Moreover, $K^-$ defined in~(\ref{eq:K_def}) can be written as $K^-=\sum_{n=0}^\infty\left(\Pi^+_-e^{B\Lambda^+}\Pi^-_+e^{B\Lambda^-}\right)^n$. These facts and the strong Markov property lead to the following corollary.
\begin{cor}\label{cor:exp_num}
It holds that
\begin{align*}
&\e\sum_{x\geq 0}\1{X(\tau^L_x)=y,J(\tau^L_x)}=\Pi^-e^{(x_0-y)\Lambda^-}K^-,&y\leq 0,\\
 &\e\sum_{x\geq
0}\1{X(\tau^L_x)=y,J(\tau^L_x)}=\Pi^+e^{(B+y-x_0)\Lambda^+}\Pi^-_+e^{B\Lambda^-}K^-,&y>0.
\end{align*}
\end{cor}

The final step of the proof is given in the following lemma.
\begin{lem}\label{lem:interchange} For any measurable non-negative function $f$ it holds a.s.~that
\begin{align*}\int_0^\infty f(X(\tau^L_x))\1{J(\tau_x^L)=j}\D x=\int_{-\infty}^\infty f(y)\sum_{x\geq 0}\1{X(\tau^L_x)=y,J(\tau^L_x)=j}\D
y.\end{align*}
\end{lem}
Intuitively this lemma states that we can interchange the
`integrals'.
\begin{proof}
Recall that $X(\tau_x^L)=U(\tau_x^L)-x$ and $U(\tau_x^L)$ is
piecewise constant. Suppose $U(\tau_x^L)=C$ for all $x\in[S,T)$
then it is immediate that
\begin{align*}&\int_{-\infty}^\infty f(y)\sum_{x\in[S,T)}\1{C-x=y,J(\tau^L_x)=j}\D
y\\
=&\int_{-\infty}^\infty f(y)\1{(C-y)\in[S,T),J(\tau^L_{C-y})=j}\D
y=\int_S^T f(C-x)\1{J(\tau_x^L)=j}\D x.\end{align*}
Summing over all such
intervals yields the statement of the lemma.
\end{proof}

\begin{proof}[Proof of Theorem~\ref{thm_localtime}]
Apply Lemma~\ref{lem:interchange} with $f(y)=e^{\a y}$ to
Corollary~\ref{cor:exp_num} to obtain
\[\int_0^\infty \e[e^{\a X(\tau_x^L)};J(\tau_x^L)]\D x=\Pi^-\int_0^\infty e^{y(\Lambda^--\a\matI)}\D ye^{x_0\Lambda^-}K^-+\Pi^+\int_0^\infty e^{y(\Lambda^++\a\matI)}\D y e^{(B-x_0)\Lambda^+}\Pi^-_+e^{B\Lambda^-}K^-.\]
Consider the MMBM $(-X(t),J(t))$ started in $B-x_0$ to find that
\[\int_0^\infty \e[e^{\a X(\tau_x^U)};J(\tau_x^U)]\D x=\Pi^+\int_0^\infty e^{y(\Lambda^++\a\matI)}\D ye^{(B-x_0)\Lambda^+}K^++\Pi^-\int_0^\infty e^{y(\Lambda^--\a\matI)}\D y e^{x_0\Lambda^-}\Pi^+_-e^{B\Lambda^+}K^+.\]
 The integrals on the right hand sides
converge if $\rho^++\a<0$ and $\rho^--\a<0$, that is
$\a\in(\rho^-,-\rho^+)$. Hence the left hand sides converge for
such $\a$. Use (\ref{lt_eq2}) to see that $k^L(\a)<0, k^U(\a)<0$
and (\ref{eq:thm1}) holds.
\end{proof}

\subsection{Stationary overflow and unused capacity}\label{sec:overflow}
Recall that $J(\tau_x^L)$ and $J(\tau_x^U)$ are irreducible
recurrent Markov chains. Denote the corresponding stationary
distributions through $\pii^L$ and $\pii^U$. Define
$\lim_{t\rightarrow\infty}L(t)/t=\kappa^L$ and similarly
$\lim_{t\rightarrow\infty}U(t)/t=\kappa^U$, which do not depend on
the initial distribution of $(X(0),J(0))$. Noting that
$\{\tau_x^L<t\}=\{x<L(t)\}$ we write
\begin{equation}\label{eq:stat_overflow}\kappa^L\pi^L_i=\lim_{t\rightarrow
\infty}\frac{1}{t}\int_0^{L(t)}\1{J(\tau_x^L)=i}\D
x=\lim_{t\rightarrow
\infty}\frac{1}{t}\int_0^\infty\1{J(\tau_x^L)=i,\tau_x^L<t}\D
x\text{ a.s.},\end{equation} i.e., one can interpret the vectors
$\kappa^L\pii^L$ and $\kappa^U\pii^U$ as the mean unused capacity
and the mean overflow in a unit of time in stationarity.

These quantities were identified in~\cite{lambda} using martingale
calculations and the theory of generalized Jordan chains.
\begin{prop}
\begin{equation}\label{eq:corr}
(\kappa^L\pii^L,\kappa^U\pii^U)\left(\begin{array}{cc}
  \matI & -\Pi^+_-e^{B\Lmb^+} \\
  -\Pi_+^-e^{B\Lmb^-} & \matI
\end{array}\right)=
\left\{\begin{array}{ll}
\kappa(\zero,\pii^+),&\hbox{ if }\kappa> 0 \\
\zero,&\hbox{ if }\kappa=0 \\
-\kappa(\pii^-,\zero),&\hbox{ if }\kappa< 0
\end{array}\right.,
\end{equation}
where $q=0$ and $\pii^\pm$ is the stationary distribution
associated with $\Lambda^\pm$ whenever $\pm\kappa>0$.
\end{prop}
This result is Corollary 5.1 of~\cite{lambda}. Note that
(\ref{eq:corr}) uniquely identifies $\kappa^L\pii^L$ and
$\kappa^U\pii^U$ if $\kappa\neq 0$. If, however, $\kappa=0$, then
an additional equation is required, see~\cite{lambda} for details.
It should be mentioned, that stationary overflow and unused
capacity for an arbitrary MAP are considered in~\cite{loss_rates}.
The results of this paper, however, depend on the stationary
distribution, which has no explicit solution unless the MAP can be
reduced to an MMBM. Moreover, a restrictive assumption about the
number of roots of a certain equation is made, which can fail even
in the case of an MMBM.

 In the following we
retrieve the above result by a simple argument using
Theorem~\ref{thm_localtime}. Observe that
$-M^L(0)=\int_0^\infty\p(\tau_x^L<e_q,J(\tau_x^L))\D x$ and hence
by~(\ref{eq:stat_overflow}) we get $-\lim_{q\downarrow
0}q\left[M^L(0)\right]_{ij}=\kappa^L\pi^L_j$ for any $i\in E,j\in
E^-$. A similar identity holds true if one considers $\tau^\pm_x$.
Firstly, note that
$-\Pi^+(\Lambda^+)^{-1}=\int_0^\infty\p(\tau_x^+<e_q,J(\tau_x^+))\D
x$ and $\{\tau_x^+<e_q\}=\{x<\overline X(e_q)\}$, where $\overline
X(t)=\sup\{X(s):0\leq s\leq t\}$. Secondly, it is well-known that
$\lim_{t\rightarrow\infty}\overline X(t)/t=\kappa\1{\kappa>0}$.
Finally, consider (\ref{eq:main}) with $q>0$ and $\a=0$, multiply
both sides by $q$ and let $q\downarrow 0$ to see that every row of
this equation reduces to~(\ref{eq:corr}).

\subsection{Special cases}\label{sec:special_cases}
An important special case arises when $E^-=E^+=E$, that is there
are no states when the process evolves deterministically (a
Brownian component is always present). In this special case
$\Pi^-=\Pi^+=\matI$ and hence (\ref{eq:thm1}) reduces to
\[M^L(\a)=\left((\Lambda^--\a\matI)^{-1}e^{x_0\Lambda^-}+(\Lambda^++\a\matI)^{-1}e^{(B-x_0)\Lambda^+}e^{B\Lambda^-}\right)\left(\matI-e^{B\Lambda^+}e^{B\Lambda^-}\right)^{-1}.\]
Letting $x_0=0$ and taking the inverse we obtain
\begin{equation}\label{williams_multi_ext}F^L(\a)=\left(e^{-B\Lambda^-}-e^{B\Lambda^+}\right)\left((\Lambda^--\a\matI)^{-1}e^{-B\Lambda^-}+(\Lambda^++\a\matI)^{-1}e^{B\Lambda^+}\right)^{-1}.\end{equation}
Moreover, observing that $X(\tau_0^L)=U(\tau^L_0)$ we find from
(\ref{eq:thm1}) that
\begin{align}\label{williams_multi}&\e_{x_0=B}[e^{\a
U(\tau_0^L)-q\tau_0^L};J(\tau_0^L)]\\&\nonumber =\left((\Lambda^--\a\matI)^{-1}+(\Lambda^++\a\matI)^{-1}\right)\left((\Lambda^--\a\matI)^{-1}e^{-B\Lambda^-}+(\Lambda^++\a\matI)^{-1}e^{B\Lambda^+}\right)^{-1}\\
&\nonumber
=\left(\Lambda^--\a\matI\right)^{-1}\left(\Lambda^++\Lambda^-\right)\left(e^{-B\Lambda^-}(\Lambda^++\a\matI)+(\Lambda^--\a\matI)e^{B\Lambda^+}\right)^{-1}\left(\Lambda^--\a\matI\right),\end{align}
where in the last step we used the fact that
$(\Lambda^++\a\matI)^{-1}$ and $e^{B\Lambda^+}$ commute. A number
of other useful transforms can be found using similar algebraic
manipulations.

Next, we restrict ourselves to the case of a single state, that
is, we consider a Brownian motion ($\sigma^2,\mu$). Without real
loss of generality it is assumed that $\sigma^2=1$. According to
(\ref{eq:lambda}), $\lambda=\Lambda^\pm$ is a solution of
$1/2\lambda^2\mp\mu\lambda-q=0$. Moreover, $\Lambda^\pm$ is
negative unless $q=0$ and $\pm\mu\geq 0$, in which case it is 0.
Thus $\Lambda^+=\mu-\gamma$ and $\Lambda^-=-\mu-\gamma$, where
$\gamma=\sqrt{\mu^2+2q}$. Now the right side of
(\ref{williams_multi}) reduces to
$-2\gamma/[e^{B(\mu+\gamma)}(\mu+\a-\gamma)-e^{B(\mu-\gamma)}(\mu+\a+\gamma)]$
and so
\begin{equation}\label{eq:analytic_cont}\e_B[e^{\a
U(\tau_0^L)-q\tau_0^L}]=\frac{e^{-B\mu}}{\cosh(B\gamma)-\frac{\mu+\a}{\gamma}\sinh(B\gamma)},\end{equation}
where $\a\in(-\mu-\gamma,-\mu+\gamma)$ according to
Theorem~\ref{thm_localtime}. Fix $q>0$ for a moment, so that
$\gamma>0$. Multiply both sides of the equation by the denominator
in the right hand side and observe that the Laplace transform is
analytic in $\Re(\a)<-\mu+\gamma$. Thus the latter equality holds
in this domain. This shows that~(\ref{eq:analytic_cont}) holds for
all $q>0,\a\leq 0$. By symmetry we obtain
\begin{align*}&\e_0[e^{-\a
L(\tau_0^U)-q\tau_0^U}]=\frac{e^{B\mu}}{\cosh(B\gamma)+\frac{\mu+\a}{\gamma}\sinh(B\gamma)},&q>0,\a\geq
0,
\end{align*} which is equation (7) in~\cite{williams_RBM}. Taking
limits on both sides we show that this equation holds true also
for $q=0$, unless $\mu=0$, in which case $\e e^{-\a
L(\tau_0^U)}=1/(1+\a B)$ by the L'H\^{o}pital's rule. It is noted
that in~\cite{williams_RBM} a very different approach is used. It
uses stochastic integration and relies on a sophisticated guess of
the right form of a certain function. Our approach, however, is
direct and is based on simple probabilistic arguments.

Let us conclude by making some additional comments about the
Brownian motion with two reflecting barriers. Firstly, the strong
Markov property of $X(t)$ implies that the nondecreasing piecewise
constant process $U(\tau^L_x)$ has memory-less jumps and
inter-arrival times, implying that it is a Poisson process with
exponential jumps. Let us confirm this and find the corresponding
rates. Note that (\ref{williams_multi_ext}) can be rewritten as
\[F^L(\a,q)=\frac{2(\frac{1}{2}\a^2+\mu\a-q)}{\gamma\coth(B\gamma)-(\mu+\a)}.\]
Hence if $\mu\neq 0$, then
\[\log\e e^{\a U(\tau^L_1)}=F^L(\a)+\a=\frac{-\a\mu(\coth(\mu B)+1)}{\a-\mu(\coth(\mu B)-1)},\]
which implies that $U(\tau^L_x)$ jumps with rate $\mu(\coth(\mu
B)+1)=2\mu/(1-e^{-2\mu B})$, and the jumps are exponential of rate
$\mu(\coth(\mu B)-1)=2\mu/(e^{2\mu B}-1)$. If $\mu=0$ then $\log\e
e^{\a U(\tau^L_1)}=\a/(1-\a B)$, that is, both rates become $1/B$.

\section*{Acknowledgements}
I would like to thank Onno Boxma and Michel Mandjes for their comments and suggestions.
I am also grateful to Yoni Nazarathy, through whom I became aware of the relevant work~\cite{williams_RBM} by Ruth Williams.
\bibliography{MMBM_2sided}

\begin{thebibliography}{10}

\bibitem{ahn_badescu_ramaswami}
{\sc Ahn, S., Badescu, A.~L. and Ramaswami, V.} (2007).
\newblock Time dependent analysis of finite buffer fluid flows and risk models
  with a dividend barrier.
\newblock {\em Queueing Systems. Theory and Applications\/} {\bf 55,} 207--222.

\bibitem{anick_mitra_sondhi}
{\sc Anick, D., Mitra, D. and Sondhi, M.~M.} (1982).
\newblock Stochastic theory of a data-handling system with multiple sources.
\newblock {\em The Bell System Technical Journal\/} {\bf 61,} 1871--1894.

\bibitem{asmussen_fluid_flow}
{\sc Asmussen, S.} (1995).
\newblock Stationary distributions for fluid flow models with or without
  {B}rownian noise.
\newblock {\em Communications in Statistics. Stochastic Models\/} {\bf 11,}
  21--49.

\bibitem{asmussen:apq}
{\sc Asmussen, S.} (2003).
\newblock {\em Applied Probability and Queues} 2nd~ed.
\newblock Applications of Mathematics. Springer-Verlag New York, Inc.

\bibitem{loss_rates}
{\sc Asmussen, S. and Pihlsg{\aa}rd, M.} (2007).
\newblock Loss rates for {L}\'evy processes with two reflecting barriers.
\newblock {\em Mathematics of Operations Research\/} {\bf 32,} 308--321.

\bibitem{breuer}
{\sc Breuer, L.} (2008).
\newblock First passage times for {M}arkov additive processes with positive
  jumps of phase type.
\newblock {\em Journal of Applied Probability\/} {\bf 45,} 779--799.

\bibitem{lambda}
{\sc D'Auria, B., Ivanovs, J., Kella, O. and Mandjes, M.}
\newblock First passage process of a {M}arkov additive process, with
  applications to reflection problems.
\newblock (submitted).

\bibitem{time_rev}
{\sc Ivanovs, J. and Mandjes, M.} (2010).
\newblock First passage of time-reversible spectrally negative {M}arkov
  additive processes.
\newblock {\em Operations Research Letters\/} {\bf 38,} 77 -- 81.

\bibitem{karandikar_kulkarni}
{\sc Karandikar, R.~L. and Kulkarni, V.~G.} (1995).
\newblock Second--order fluid flow models: Reflected {B}rownian motion in a
  random environment.
\newblock {\em Operations Research\/} {\bf 43,} 77--88.

\bibitem{kosten}
{\sc Kosten, L.} (1974/75).
\newblock Stochastic theory of a multi-entry buffer. {I}.
\newblock {\em Delft Progress Report\/} {\bf 1,} 10--18.

\bibitem{ramanan_skorohod0a}
{\sc Kruk, L., Lehoczky, J., Ramanan, K. and Shreve, S.} (2007).
\newblock An explicit formula for the {S}korokhod map on {$[0,a]$}.
\newblock {\em The Annals of Probability\/} {\bf 35,} 1740--1768.

\bibitem{fluid_survey}
{\sc Kulkarni, V.~G.} (1997).
\newblock Fluid models for single buffer systems.
\newblock In {\em Frontiers in queueing}.
\newblock Probability and Stochastics Series. CRC pp.~321--338.

\bibitem{lindley59}
{\sc Lindley, D.} (1959).
\newblock Discussion on {M}r. {W}insten's paper.
\newblock {\em Journal of the Royal Statistical Society\/} {\bf B21,} 22--23.

\bibitem{pistorius}
{\sc Pistorius, M.} (2006).
\newblock On maxima and ladder processes for a dense class of {L}\'evy process.
\newblock {\em Journal of Applied Probability\/} {\bf 43,} 208--220.

\bibitem{prabhu_zhu}
{\sc Prabhu, N.~U. and Zhu, Y.} (1989).
\newblock {M}arkov-modulated queueing systems.
\newblock {\em Queueing Systems. Theory and Applications\/} {\bf 5,} 215--245.

\bibitem{rogers}
{\sc Rogers, L. C.~G.} (1994).
\newblock Fluid models in queueing theory and {W}iener-{H}opf factorization of
  {M}arkov chains.
\newblock {\em Annals of Applied Probability\/} {\bf 4,} 390--413.

\bibitem{williams_RBM}
{\sc Williams, R.~J.} (1992).
\newblock Asymptotic variance parameters for the boundary local times of
  reflected {B}rownian motion on a compact interval.
\newblock {\em Journal of Applied Probability\/} {\bf 29,} 996--1002.

\end{thebibliography}

\end{document}